\documentclass{amsart}
\begin{document}
\def \tA{\tilde{A}}
\def \tB{\tilde{B}}
\def \tC{\tilde{C}}
\def \alph{A}
\def \bet{B}
\def \bsigma{\bar{\sigma}}
\def \y{^{\infty}}
\def \Ra{\Rightarrow}
\def \uBS{\overline{BS}}
\def \lBS{\underline{BS}}
\def \lB{\underline{B}}
\def \<{\langle}
\def \>{\rangle}
\def \hL{\hat{L}}
\def \cU{\mathcal U}
\def \cF{\mathcal F}
\def \S{\Sigma}
\def \st{\stackrel}
\def \sd{Spec_{\d}\ }
\def \pd{Proj_{\d}\ }
\def \s{\sigma_2}
\def \i{\sigma_1}
\def \bs{\bigskip}
\def \cD{\mathcal D}
\def \cC{\mathcal C}
\def \cT{\mathcal T}
\def \cK{\mathcal K}
\def \cX{\mathcal X}
\def \sX{X_{set}}
\def \cY{\mathcal Y}
\def \cS{X}
\def \cR{\mathcal R}
\def \cE{\mathcal E}
\def \tcE{\tilde{\mathcal E}}
\def \cP{\mathcal P}
\def \cA{\mathcal A}
\def \cV{\mathcal V}
\def \cM{\mathcal M}
\def \cN{\mathcal N}
\def \tcM{\tilde{\mathcal M}}
\def \caS{\mathcal S}
\def \cG{\mathcal G}
\def \cB{\mathcal B}
\def \tG{\tilde{G}}
\def \cF{\mathcal F}
\def \h{\hat{\ }}
\def \hp{\hat{\ }}
\def \tS{\tilde{S}}
\def \tP{\tilde{P}}
\def \tA{\tilde{A}}
\def \tX{\tilde{X}}
\def \tcS{\tilde{X}}
\def \tT{\tilde{T}}
\def \tE{\tilde{E}}
\def \tV{\tilde{V}}
\def \tC{\tilde{C}}
\def \tI{\tilde{I}}
\def \tU{\tilde{U}}
\def \tG{\tilde{G}}
\def \tu{\tilde{u}}
\def \chu{\check{u}}
\def \tx{\tilde{x}}
\def \tL{\tilde{L}}
\def \tY{\tilde{Y}}
\def \d{\delta}
\def \e{\epsilon}
\def \bZ{{\bf Z}}
\def \bV{{\bf V}}
\def \bF{{\bf F}}
\def \bE{{\bf E}}
\def \bC{{\bf C}}
\def \bO{{\bf O}}
\def \bR{{\bf R}}
\def \bA{{\bf A}}
\def \bB{{\bf B}}
\def \cO{\mathcal O}
\def \ra{\rightarrow}
\def \bx{{\bf x}}
\def \f{{\bf f}}
\def \bX{{\bf X}}
\def \bH{{\bf H}}
\def \bS{{\bf S}}
\def \bF{{\bf F}}
\def \bN{{\bf N}}
\def \bK{{\bf K}}
\def \bE{{\bf E}}
\def \bB{{\bf B}}
\def \bQ{{\bf Q}}
\def \bd{{\bf d}}
\def \bY{{\bf Y}}
\def \bU{{\bf U}}
\def \bL{{\bf L}}
\def \bQ{{\bf Q}}
\def \bP{{\bf P}}
\def \bR{{\bf R}}
\def \bC{{\bf C}}
\def \bD{{\bf D}}
\def \bM{{\bf M}}
\def \bP{{\bf P}}
\def \xtoleqr{x^{(\leq r)}}
\def \hU{\hat{U}}

\newtheorem{THM}{{\!}}[section]
\newtheorem{THMX}{{\!}}
\renewcommand{\theTHMX}{} 
\newtheorem{thm}{Theorem}[section]
\newtheorem{cor}[thm]{Corollary}
\newtheorem{conj}[thm]{Conjecture}
\newtheorem{lem}[thm]{Lemma}
\newtheorem{question}[thm]{Question}
\newtheorem{prop}[thm]{Proposition}
\theoremstyle{definition}
\newtheorem{defn}[thm]{Definition}
\theoremstyle{remark}
\newtheorem{rem}[thm]{Remark}
\newtheorem{example}[thm]{\bf Examples}
\numberwithin{equation}{section}
\title{Complex dynamics and invariant forms mod $p$}
\author{Alexandru Buium}
\address{University of New Mexico \\ Albuquerque, NM 87131}
\email{buium@math.unm.edu}
\maketitle
\begin{abstract}
Complex dynamical systems on the Riemann sphere do not possess
``invariant
forms''. However there exist non-trivial examples
of dynamical systems, defined over number fields,
satisfying the property that their reduction modulo $\wp$
 possesses ``invariant forms'' for all but finitely many places
$\wp$. The paper completely characterizes
the dynamical systems possessing the latter property.
\end{abstract}

\section{Motivation and statement of the Theorem}

Let $k$ be an algebraically closed field; in our applications
$k$ will be either the complex field $\bC$ or the algebraic closure
 of a finite field. 
Let $k(t)$
be the field of rational functions in the variable $t$
and let $\sigma(t) \in k(t)$ be a non-constant
rational function.
 By an {\it invariant form of weight $\nu \in \bZ$} 
for $\sigma(t)$ we mean a rational function
   $f(t) \in k(t)$
such that
\begin{equation}
\label{viinnt}
f(\sigma(t)) = \left( \frac{d \sigma}{dt} (t) \right)^{-\nu}
\cdot
 f(t).
\end{equation}
If we consider the  $\nu-$tuple differential form
  $\omega:=f(t) \cdot (dt)^{\nu}$ then 
Equation
\ref{viinnt} can be written  as $\sigma^* \omega=\omega$
which is actually what justifies our terminology. We shall sometimes
refer to $\omega$ itself as being an {\it invariant form of weight
  $\nu \in \bZ$}.

Our first remark is that if $k=\bC$ and $\sigma$, viewed as a self
map $\bP^1
\ra \bP^1$ of the complex projective line,
 has degree $\geq 2$ then there are no non-zero invariant
forms for $\sigma$ of non-zero weight. Indeed if Equation \ref{viinnt}
holds for some $f \neq 0$ and $\nu \neq 0$
 then the same equation holds with $\sigma$ replaced
by the $n-$th iterate $\sigma^n$ for any $n$. Now this equation
for the iterates implies that any finite non-parabolic periodic point of
$\sigma$ 
is either a zero or a pole of $f$; cf. \cite{Milnor}, p. 99 for the
definition
of parabolic periodic points. On the other hand
by \cite{Milnor}, pp. 47 and 143,
there are  
infinitely many non-parabolic periodic points, a contradiction.

Although invariant forms don't exist over the complex numbers
there exist, nevertheless, interesting examples 
of complex rational functions with coefficients in number fields
whose reduction mod (almost all)
primes
admit invariant forms.
 To explain this
let $F \subset \bC$ be a number field
(always assumed of finite degree over the rationals)
and let $\wp$ be a place of $F$ (always assumed finite). 
Let $\cO_{\wp}$ denote the valuation ring of $\wp$, let 
$\kappa_{\wp}$ be the residue field of $\cO_{\wp}$
and let 
$k_{\wp}=\kappa^a_{\wp}$ be an algebraic closure of
$\kappa_{\wp}$.
Now if  
$\sigma(t) \in \bC(t)$  has coefficients in  $F$
then  for
 all except finitely many places  $\wp$ of $F$ we may consider
the rational function $\bsigma_{\wp}(t) \in \kappa_{\wp}(t) \subset
k_{\wp}(t)$ obtained
by writing $\sigma(t)=P(t)/Q(t)$, 
$P(t),Q(t) \in \cO_{\wp}[t]$, with $Q(t)$
 primitive, and then reducing the coefficients
of $P(t)$ and $Q(t)$ modulo the maximal ideal of $\cO_{\wp}$. 
  Below is a list of  examples of $\sigma(t)$'s such that
  $\bsigma_{\wp}(t)$
has non-zero invariant forms of non-zero weight for all but finitely
  many places $\wp$ of $F$. Our main result will then state that
the examples below are {\it all} possible examples.

\begin{example}
\label{scarpin}
1) {\it Multiplicative functions}. Let $F={\bf Q}$ and
let $\sigma(t)=t^{\pm d}$ where $d$ is a positive integer. Then for all 
prime integers $p$ not dividing $d$ the rational function
$\bsigma_p(t) \in \bF_p^a[t]$ 
 possesses a non-zero
  invariant form  of weight $p-1$,
\[\omega_p=\left(\frac{dt}{t} \right)^{p-1}.\]

2) {\it Chebyshev polynomials}.
Let $F={\bf Q}$ and let 
$Cheb_d(t) \in \bC[t]$ be the Chebyshev polynomial of degree $d$, i.e.
the unique polynomial in $\bC[t]$ such that
$Cheb_d(t+t^{-1})=t^d+t^{-d}$. Clearly $Cheb_d$ has integer coefficients.
Let $\sigma(t)=\pm Cheb_d(t)$.
Then for all odd prime integers $p$ not dividing $d$ the polynomial
$\bsigma_p(t) \in \bF_p^a[t]$ 
 possesses a non-zero
  invariant form  of weight $p-1$, 
\[\omega_p=\frac{(dt)^{p-1}}{(t^2-4)^{(p-1)/2}}.\]

3) {\it Latt\`{e}s functions} \cite{Lattes},
\cite{Milnor}. 
A  rational
 function $\sigma(t) \in \bC(t)$,
$\sigma:\bP^1 \ra \bP^1$, will be called a {\it 
  Latt\`{e}s
function}    if
 $\sigma$ is obtained as follows.
One starts with  an elliptic curve $E$
over $\bC$ with affine plane equation $y^2=x^3+ax+b$.
One then considers an 
algebraic group endomorphism
$\tau_0:E \ra E$, a non-trivial
 algebraic group
automorphism $\gamma_0:E \ra E$, and a fixed point
$P_0$ of $\gamma_0$.
Note that $\gamma_0$ has order $2,4,3$ or $6$
and, correspondingly, one has:
$\gamma_0(x,y)=(x,-y)$ or $b=0$, $\gamma_0(x,y)=(-x,iy)$,
or $a=0$, $\gamma_0(x,y)=(\zeta_3x,y)$, or
$a=0$, $\gamma_0(x,y)=(\zeta_3x,-y)$.
 One considers the morphism
$\tau:E \ra E$, $\tau(P)=\tau_0(P)+P_0$ and one considers
the  map  
 $\sigma:E/\Gamma  \ra E/\Gamma $ induced by $\tau$,
where $\Gamma:=\<\gamma_0\>$ and
$E/\Gamma$ is identified with $\bP^1=Proj\ \bC[x_0,x_1]$ via the 
isomorphism that sends $t=x_1/x_0$ into 
$x,x^2,y$, or $y^2$ according as $ord(\gamma_0)$ is $2,4,3$ or $6$.
We denote by $Lat_d$ any Latt\`{e}s function of degree $d$.
The original examples of Latt\`{e}s, cf. \cite{Milnor},
are given by the  above construction
with  $ord(\gamma_0)=2$.

We claim that if $\sigma(t)$ is a 
 Latt\`{e}s function with  $E$
 defined over a number field $F$ then $\bsigma_{\wp}(t)$ has a
  non-zero invariant form of non-zero weight for all but finitely many
  places $\wp$ of $F$. Indeed if $\omega_E$ is a non-zero global $1-$
  form on $E$
defined over $F$ and $e=ord(\gamma_0)$ then $\omega_E^e=\pi^* \omega$
where $\pi:E \ra E/\<\gamma_0\>$ is the canonical projection and
 $\omega
 \in
F(t) \cdot (dt)^e$. Now 
\[\pi^* \sigma^* \omega=\tau^* \pi^* \omega=(\tau^* \omega_E)^e=\lambda
 \omega_E^e=
\pi^* (\lambda \omega)\]
for some $\lambda \in F^{\times}$. Hence $\sigma^* \omega=\lambda
\omega$
 hence 
\[\bsigma_{\wp}^* \bar{\omega}_{\wp}^{q_{\wp}-1}=
\bar{\omega}_{\wp}^{q_{\wp}-1}\]
for all but finitely many $\wp$'s, where 
$\bar{\omega}_{\wp}$ is the reduction mod $\wp$ of $\omega$ and
$q_{\wp}$ is the size of the
residue field  $\kappa_{\wp}$.

4) {\it Flat functions}.
It is convenient to introduce terminology
that ``puts together'' all the above Examples. First
let us say that two rational functions $\sigma_1, \sigma_2 \in \bC(t)$
are {\it conjugate} if there exists $\varphi \in \bC(t)$ of degree $1$
such that $\sigma_2=\varphi \circ \sigma_1 \circ \varphi^{-1}$; if this is
the case
write $\sigma_1 \sim \sigma_2$. Say that $\sigma(t) \in \bC(t)$ is
{\it flat} if either $\sigma(t) \sim t^{\pm d}$ or
$\sigma(t) \sim \pm Cheb_d(t)$ or $\sigma(t) \sim Lat_d$.
(This is an ad hoc  version of  terminology from
\cite{Milnor2}.)
Flat rational maps naturally appear
 in a number of  contexts; cf Remark \ref{cucurigu}.
In the present paper
our interest in flat functions comes from the fact that
 if $\sigma(t) \in \bC(t)$ is flat
and has coefficients in a number field $F$ then
$\bsigma_{\wp}(t)$ has a
  non-zero invariant form of non-zero weight for all but finitely many
  places $\wp$ of $F$. To check this
 we may write $\sigma=\varphi \circ \sigma_1
 \circ  \varphi^{-1}$
where  $\sigma_1 = t^{\pm d}$ or $\sigma_1=\pm Cheb_d(t)$ or $\sigma_1=Lat_d$
and $\varphi \in \bC(t)$, $deg(\varphi)=1$,
cf. Examples $1,2,3$ above. A standard specialization argument shows that
we may replace $\sigma_1$ and $\varphi$ by maps defined over a finite
extension $F'$ of $F$, where, in case $\sigma_1=Lat_d$,
 the elliptic curve to
which $\sigma_1$ is attached  
is also defined over $F'$. Then we may conclude
by the discussion of the Examples  $1,2,3$.
\end{example}

The aim of this note is to prove  that, conversely, we have the following:

\begin{thm}
\label{cococo}
Assume $\sigma(t) \in \bC(t)$ 
has coefficients in a number field $F$ and has degree $d \geq 2$.
Assume that,
for infinitely many places $\wp$ of $F$,
 $\bsigma_{\wp}(t) \in k_{\wp}(t)$ admits a non-zero invariant
form of non-zero weight.
 Then $\sigma$ is flat.
\end{thm}

We stress the fact that the form and the weight in the above statement
depend a priori on $\wp$.

The proof of the Theorem  will be presented in Section 2. 
The idea is to 
show that the existence of invariant forms for $\bsigma_{\wp}$
 implies that  $\bsigma_{\wp}$
have a ``very special type of ramification''.
 This will imply that
 $\sigma$ itself has the same type of ramification.
Then
one concludes 
by  using the topological
characterization
of postcritically finite rational maps with non-hyperbolic orbifold
 \cite{DH}; the latter is, itself, an easy application of Thurston's
orbifold theory.

In the rest of Section 1 we make some Remarks and raise some questions;
 in particular
we indicate the connection between our Theorem and 
results in  \cite{shimura}, \cite{dyn},
\cite{Ihara}.

\begin{rem}
Our Theorem suggests, more generally, the problem of characterizing
all pairs  $(\sigma_1,\sigma_2)$ of non-constant
morphisms $\sigma_1,\sigma_2:\tY \ra Y$
of smooth projective curves over a number field $F$ satisfying the
property
that:

\bigskip

(*) {\it  For all except finitely 
 many places $\wp$ of $F$, there exists an integer
$\nu_{\wp} \neq 0$ and a non-zero $\nu_{\wp}-$tuple form
$\omega_{\wp} \in \Omega^{\otimes \nu_{\wp}}_{k_{\wp}(Y)/k_{\wp}}$
with the property that}
$\bsigma_{1,\wp}^* \omega_{\wp}=\bsigma_{2,\wp}^* \omega_{\wp}$.

\bigskip

Here $k_{\wp}(Y)$ is the field of rational functions on the reduction
of $Y$ mod $\wp$ and $\bsigma_{i,\wp}$ are the induced
morphisms mod $\wp$ (with respect to some model of $\sigma_i$
over a ring of $S-$integers of $F$, $S$ a finite set of places of $F$). 
The only interesting case of this problem is, of course, that in which
the smallest equivalence relation in $Y \times Y$ containing the image
of
\[\sigma_1 \times \sigma_2:\tY \ra Y \times Y\]
is Zariski dense in $Y \times Y$.
The Examples \ref{scarpin}
fit
 into the above scheme with $Y=\tY=\bP^1$ and
$\sigma_1=id$. If one allows, however, both $\sigma_i$'s to be
different from the identity
then there are interesting additional examples of pairs $(\sigma_1,\sigma_2)$
satisfying property (*). A series of 
such examples  originates in the work of
\cite{Ihara}: the curves $Y$ and $\tY$ are in this case modular or
 Shimura curves
and the forms $\omega_{\wp}$ correspond to ``Hasse invariants at
$\wp$''
viewed as appropriate 
modular forms mod $\wp$. 
Other, more elementary, examples of pairs
$(\sigma_1,\sigma_2)$
satisfying property (*) can be produced by taking, for instance,
$Y=\tY=\bP^1$ and $\sigma_2=\sigma_1 \circ \tau$ where $\sigma_1:
\bP^1 \ra \bP^1$ is a
Galois cover and $\tau$ is an automorphism of $\bP^1$.
 One can wonder if these
examples, together with the Examples \ref{scarpin}
 are (essentially) the only examples for which (*)
is satisfied. Cf. also \cite{CU} for a measure theoretic analogue.
\end{rem}

\begin{rem}
Let us say that a 
 $\nu-$tuple form $\omega \in k(t)(dt)^{\nu}$ of weight $\nu$
is  {\it semi-invariant} 
for  $\sigma$ if $\sigma^* \omega=\lambda
\omega$ for some $\lambda \in k^{\times}$. 
One can show that $\sigma(t) \in \bC(t)$ is flat
if and only if there
exists a non-zero form $\omega \in \bC(t) \cdot
(dt)^{\nu}$ of weight 
$\nu \neq 0$ which is {\it semi-invariant} for $\sigma$.
 (For the ``only if''
part one can  use  the arguments in the discussion of the Examples
\ref{scarpin};  the ``if'' part follows, for instance,
from Lemmas \ref{sase}, \ref{sapte} and Proposition \ref{DHTHM}
below.)
So our Theorem implies, in particular, the following statement:

\medskip

{\it (**)
If $\sigma \in \bC(t)$ has coefficients in a number field and 
$\bar{\sigma}_{\wp}(t)$ possesses a 
non-zero invariant form of non-zero weight
for  infinitely many places $\wp$ then $\sigma(t)$ itself
possesses a
non-zero semi-invariant
form of non-zero weight.}

\medskip

 On the other hand  (**) fails
in the context of correspondences! Indeed if
$\sigma_1,\sigma_2:\tY \ra Y$ define a Hecke correspondence between
modular curves over the complex numbers
it generally  happens that condition (*) in the
preceding
Remark holds but, nevertheless, there is no non-zero
{\it semi-invariant}
form of non-zero weight
for this correspondence in characteristic zero i.e. no 
non-zero $\nu-$tuple form 
$\omega \in \Omega^{\otimes \nu}_{\bC(Y)/\bC}$, $\nu \neq 0$,
with the property that
$\sigma_{2}^* \omega= \lambda \cdot \sigma_{1}^* \omega$
for some $\lambda \in \bC^{\times}$.
\end{rem}

\begin{rem}
\label{cucurigu}
Flat rational maps naturally appear
 in a number of (a priori unrelated) contexts such as:
dynamical systems with smooth Julia sets \cite{Milnor}, pp
67-70,
orbifold theory
 \cite{DH}, \cite{Milnor}, density conjecture for hyperbolic rational maps
\cite{Mane}, isospectral deformations of rational maps
\cite{McMullen}, permutable rational maps \cite{Ritti}, measure theory
related to rational maps \cite{Zdunick}, and
Galois theory of function fields related to ``Shur's conjecture''
\cite{fried}, \cite{guralnick}.
Cf. also \cite{Milnor2} for an overview of some of these topics.

The result of the present paper is closely related (and can be applied
to)
the
theory developed in \cite{dyn}, \cite{shimura}, \cite{buc}.
That theory implies the existence of 
 what one can call ``invariant $\d-${\it
 forms}''
 for flat functions
\cite{dyn} and for
 Hecke correspondences on modular and Shimura curves \cite{shimura}.
These invariant $\d-$forms are objects in  characteristic zero
but transcending usual algebraic geometry (because they involve
not only the coordinates of the points but also their iterated ``Fermat
quotients''
up to a certain ``order''). It turns out that the presence of
invariant 
$\d-$forms
implies the presence of invariant forms mod $p$ 
(in the sense of the present paper);
this allows one to prove ``converse theorems'' along the lines
of the ``main questions'' raised in \cite{dyn}. Explaining this 
application would require reviewing
the  $\d-$geometric context of the above papers.
 This would go beyond the scope of the
present paper and will be discussed elsewhere \cite{buc}. 
\end{rem}

{\it Acknowledgments.} The author would like to thank Y. Ihara and
J. H. Silverman for their comments and suggestions on a preliminary
version of this paper. The preliminary version was based on an elementary
(but somewhat tedious) analysis of the poles and zeroes of invariant
forms. The more conceptual proof in the present version of the paper
was suggested to the author by Y. Ihara.

\section{Proof of the Theorem}

For the proof we need some preliminaries.
Let $k$ be  an algebraically closed 
 field  of  characteristic $p \geq 0$.
For later purposes we will assume $p \neq 2,3$.
We identify the projective line $\bP^1$ over $k$
with its set
$k \cup \{\infty\}$
 of $k-$points.
 Let  $\sigma \in k(t)$ be a
 non-constant
rational
 function which we view as a rational function
$\sigma:\bP^1
\ra \bP^1$.
We denote by $e_{\sigma}:\bP^1 \ra \bZ_{>0}$ the ramification
index function: for $B \in \bP^1$, $e_{\sigma}(B)$ is the valuation
at $B$, $v_B(\sigma^* t_A)$,  of $\sigma^* t_A$ where $t_A$ is a parameter 
of the local ring of $\bP^1$ 
at $A=\sigma(B)$. Recall that $\sigma$ is {\it tamely ramified}
if $e_{\sigma}(B)$ is not divisible by $p$  for all $B$.
Recall \cite{DH}
that one defines the {\it critical locus} $\Omega_{\sigma}$
and {\it
 postcritical locus} $P_{\sigma}$ of $\sigma$ as
\[\Omega_{\sigma}:=\{B \in \bP^1\ |\ 
e_{\sigma}(B) >1\}\ \ \text{and} \ \ P_{\sigma}:=\bigcup_{n \geq 1}
 \sigma^n(\Omega_{\sigma})\]
respectively. 
For any integer $n \geq 1$ we have $P_{\sigma}=P_{\sigma^n}$.
If $\sigma$ is separable then $\Omega_{\sigma}$ is finite;
if in addition
 $k$ is the algebraic closure of a finite field
then $P_{\sigma}$ is also finite.
(However, if $char(k)=0$, $P_{\sigma}$ is generally infinite.)

The function $e_{\sigma}$ and the critical locus $\Omega_{\sigma}$
make sense, of course, in the more general situation when $\sigma$ is
a morphism of non-singular curves over $k$.

Assume now  $\sigma(t) \in k(t)$ has
degree $\geq 2$;
then we define, following \cite{DH}, \cite{Milnor},
a function $\mu_{\sigma}:\bP^1 \ra \bZ_{>0} \cup \{\infty\}$ by the formula
\[\mu_{\sigma}(A):=lcm\ \{e_{\sigma^m}(B)\ |\ m \geq 1, \ \sigma^m(B)=A\}.\]
Here $lcm$ stands for 
{\it lowest common multiple} in the multiplicative monoid
$\bZ_{>0} \cup \{\infty\}$.
Clearly  $\mu_{\sigma}(A)=1$ for $A \not\in P_{\sigma}$
and $\mu_{\sigma}(A)>1$ for  $A \in
P_{\sigma}$.
The function $\mu_{\sigma}$
can be characterized as being the smallest among all functions
$\mu:\bP^1 \ra \bZ_{>0} \cup \{\infty\}$ 
such that $\mu(A)=1$ for $A \not\in P_{\sigma}$ and such that
$\mu(A)$ is a multiple of $\mu(B) \cdot e_{\sigma}(B)$ for each $B \in
\sigma^{-1}(A)$. 
One defines the {\it orbifold} attached to $\sigma$ as the pair
$O_{\sigma}:=(\bP^1,\mu_{\sigma})$
 and  one defines the {\it Euler characteristic}
as being the rational number
\[\chi(O_{\sigma}):=2-\sum_{A \in P_{\sigma}}
\left( 1- \frac{1}{\mu_{\sigma}(A)} \right) \in \bQ.\]
Over the complex numbers we have:

\begin{prop}
\label{DHTHM}
\cite{DH}
Assume $\sigma(t)
 \in \bC(t)$ has degree $\geq 2$,  $P_{\sigma}$ is finite, and
 $\chi(O_{\sigma})=0$. Then 
$\sigma$ is flat.
\end{prop}

\begin{proof}
This follows by combining Proposition 9.2, p. 290,
Proposition 9.3, p. 290, and
Corollary
2.4, p. 269, in \cite{DH}.
\end{proof}

In what follows we will prove a series of Lemmas. 
 Recall that a
 $\nu-$tuple form $\omega \in \Omega_{k(X)/k}^{\otimes \nu}$
on a non-singular projective curve $X$ over $k$ 
is   called {\it semi-invariant} (of weight $\nu$)
for an endomorphism $\sigma:X \ra X$ if $\sigma^* \omega=\lambda
\omega$ for some $\lambda \in k^{\times}$.

\begin{lem}
\label{unu}
Assume $\sigma(t) \in k(t)$ is tamely ramified.
 Let $B \in \bP^1$, $A:=\sigma(B)$, and let $\omega=f(dt)^{\nu}$
 be a  $\nu-$tuple form.
Then
\begin{equation}
\label{ihara1}
ord_B(\sigma^* \omega)+\nu=e_{\sigma}(B) \cdot (ord_A(\omega)+\nu).
\end{equation} 
\end{lem}

\begin{proof}
One may assume
\[\sigma^* t_A=t_B^e+... \in k[[t_B]],\ \ \ \omega=
(t_A^n+...)(dt_A)^{\nu} \in k((t_A))(dt_A)^{\nu},\]
where  $e=e_{\sigma}(B)$, $n=ord_A(\omega)$.
Then one gets
\[\sigma^* \omega=(t_B^{en}+...)(e 
t_B^{e-1}+...)^{\nu}(dt_B)^{\nu}=(e^{\nu}t_B^{en+e
  \nu -\nu}+...)(dt_B)^{\nu},\]
and we are done.
\end{proof}

\begin{lem}
\label{doi}
If $\sigma(t) \in k(t)$ is tamely ramified of degree $d \geq 2$
and admits a semi-invariant non-zero 
form of weight one then $\sigma \sim t^{\pm d}$. 
\end{lem}

Here $\sim$ means conjugation by an automorphism of $k(t)$.

\begin{proof}
Assume $\sigma^* \omega=\lambda \omega$, $\lambda \in k^{\times}$.
Let $S$ be the support of the divisor $(\omega)$. Let $B \not\in S$,
$A :=\sigma(B)$. By
Equation \ref{ihara1} we have $1=e_{\sigma}(B)(ord_A(\omega)+1)$ hence
$e_{\sigma}(B)=1$ so $ord_A(\omega)=0$ so $A \not\in S$. This shows
that
$\sigma^{-1}(S) \subset S$. It follows that the restriction
$\sigma_{|S}:S \ra S$
is a bijection  and $\sigma^*[S]=d[S]$. So there is an
integer $N \geq 1$ such that $\sigma^N_{|S}$ is the identity on $S$.
By Equation \ref{ihara1} we get
\[ord_A(\omega)+1=e_{\sigma^N}(A)(ord_A(\omega)+1)=d^N(ord_A(\omega)+1)\]
for all $A \in S$. Hence $ord_A(\omega)=-1$ for all $A \in S$.
Since the genus of $\bP^1$ is zero, $deg\ (\omega)=-2$. Hence
$\sharp S=2$ so we may assume $S=\{0, \infty\}$. We find that
$\sigma^{-1}(\{0,\infty\})=\{0,\infty\}$ which implies that $\sigma=c
t^{\pm d}$ for some $c \in k^{\times}$ hence $\sigma \sim t^{\pm d}$
 and we are done.
\end{proof}

\begin{lem}
\label{trei}
If $\sigma(t) \in k(t)$  admits a semi-invariant non-zero 
form of non-zero weight then $\sigma$ also admits 
a semi-invariant non-zero 
form of  positive weight not divisible by $p$.
\end{lem}

\begin{proof}
We may assume $p>0$.
Let $\omega=f(dt)^{\nu}$ be a semi-invariant non-zero 
form of non-zero weight $\nu$ where $\nu$ 
has minimum absolute value. We may assume $\nu >0$.
We claim that $\nu$ is not divisible by $p$. Assume it is and seek a
contradiction. 
If $df=0$ then $f=g^p$ for some $g \in k(t)$ and then $g(dt)^{\nu/p}$
is semi-invariant, a contradiction. So we may assume $df \neq 0$.
Since $\omega$ is semi-invariant and $\nu$ is divisible by $p$ it
follows that $\sigma^* f=u^p f$, where $u=\lambda^{1/p}(d
\sigma/dt)^{-\nu/p}$.
It follows that 
\[\sigma^* \left( \frac{df}{f} \right)=\frac{d(\sigma^* f)}{\sigma^*
  f}=
\frac{d(u^pf)}{u^pf}=\frac{u^p df}{u^p f}=\frac{df}{f}\]
so $df/f$ is invariant of weight $1<\nu$, a contradiction. 
\end{proof}

Let us consider the following construction. Assume $\sigma \in k(t)$
admits
a non-zero semi-invariant form $\omega=f(dt)^{\nu}$
 of positive weight $\nu$ not divisible
by $p$. We let 
\begin{equation}
\label{xomega}
\pi:X_{\omega} \ra \bP^1
\end{equation}
 be the cyclic Galois cover of non-singular projective curves
corresponding to the field extension 
\[k(t) \subset k(t,f^{1/\nu}).\]

\begin{lem}
\label{patru}
Assume $\sigma \in k(t)$
admits
a non-zero semi-invariant form $\omega$
 of positive weight $\nu$ not divisible
by $p$.
There exists a rational $1-$form $\eta$ on $X_{\omega}$ such that
$\pi^* \omega=\eta^{\nu}$. Moreover $\sigma:\bP^1 \ra \bP^1$ lifts
to an endomorphism $\tau:X_{\omega} \ra X_{\omega}$ admitting
$\eta$ as a semi-invariant form of weight one.
\end{lem}

\begin{proof}
Set $\eta:=f^{1/\nu} dt$; then clearly $\pi^* \omega=\eta^{\nu}$. 
Now extend $\sigma^*$ to a field embedding, $\tau^*$, of $k(t,f^{1/\nu})$
into an algebraic closure of $k(t)$. If $\sigma^* \omega=\lambda
\omega$ then we have
\[(\tau^* f^{1/\nu})^{\nu}=\sigma^* f=(\lambda^{1/\nu}(d
\sigma/dt)^{-1}f^{1/\nu})^{\nu}\]
for some $\nu-$th power $\lambda^{1/\nu}$ of $\lambda$, hence 
\[\tau^* f^{1/\nu} =\zeta \lambda^{1/\nu} \left(
\frac{d \sigma}{dt} \right)^{-1} f^{1/\nu},\]
for some $\nu-$th root of unity $\zeta$.
So $\tau^*$ induces an endomorphism of the field $k(t,f^{1/\nu})$ and
therefore it comes from an enodomorphism $\tau$ of
$X_{\omega}$. Clearly $\eta$ is a semi-invariant for $\tau$.
\end{proof}

If, in Lemma \ref{patru}, $\sigma$ has degree $\geq 2$ then $\tau$ has
degree $deg(\tau)=deg(\sigma) \geq 2$ so, by Hurwitz' formula
\cite{Hart}, p. 301,
 $X_{\omega}$ has  genus
\begin{equation}
\label{genus}
g(X_{\omega})=0\ \ or\ \ 1.
\end{equation}
Consider the Galois group
 \begin{equation}
\label{gamma}
\Gamma=\< \gamma \>
\end{equation}
 of the cover
\ref{xomega}.
If $g(X_{\omega})=1$ then $\gamma$ must have a fixed point which we
may take as the zero element for an algebraic
 group structure on $X_{\omega}$.
This being done, and assuming $p \neq 2,3$,
$\gamma$ has order $2,3,4,$ or $6$. The order $4$ can only occur for
$j-$invariant
$j=1728$ while the orders $3$ and $6$ can only occur if $j=0$. Cf.
\cite{Hart},
p. 321.

\begin{lem}
\label{sase}
Assume $\sigma:\bP^1 \ra \bP^1$  has degree $d \geq 2$ and admits a
non-zero semi-invariant form $\omega$ of positive 
weight not divisible by $p$. Assume either $p>d$ or $p=0$.
Assume moreover that $g(X_{\omega})=0$. Then  $\sharp P_{\sigma} \leq
3$
and $\chi(O_{\sigma})=0$.
\end{lem}

\begin{proof}
Consider the lifting $\tau:X_{\omega}=\bP^1 \ra \bP^1$
 and the form $\eta$ in Lemma \ref{patru}. Note that $\tau$ has the
 same degree $d$ as $\sigma$ so it is tamely ramified. By Lemma
 \ref{doi} we may assume
$\tau(t)= t^{\pm d}$. With $\pi$ and $\gamma$ defined in Equations
\ref{xomega} and \ref{gamma} note that the equation $\pi \circ
 \tau=\sigma \circ \pi$
implies that $\tau \circ \gamma=\gamma^s \circ \tau$ for some $s$.
We get  equalities of divisors
\[\tau^* \gamma^* 0+\tau^* \gamma^* \infty=d\gamma^{s*}0+
d\gamma^{s*}\infty,\]
hence $\{\gamma^{-1}(0),\gamma^{-1}(\infty)\}=\{0,\infty\}$ hence
$\gamma(t)=\lambda t$ or $\gamma(t)=\lambda t^{-1}$ 
for some $\lambda \in k^{\times}$.
If $\gamma(t)=\lambda t$, $n=\sharp \Gamma$, we may assume
$\pi(t)=t^n$ so $(t^{\pm d})^n=\sigma(t^n)$ so $\sigma(t)=t^{\pm d}$
hence, trivially,  $\sharp P_{\sigma}=2$ and $\chi(O_{\sigma})=0$.
If $\gamma(t)=\lambda t^{-1}$ then $\gamma^2=id$ hence we may assume
$\pi(t)=t+\lambda t^{-1}$. We must have $\tau \circ \gamma=\gamma
\circ \tau$ from which we get $\lambda^{\pm d}=\lambda$. Hence we have a
functional
equation
\[t^{\pm d}+\lambda^{\pm d} t^{\mp d}=\sigma(t+\lambda t^{-1}).\]
Taking derivatives in this equation one trivially sees that 
$\sharp P_{\sigma}=3$ and
$\chi(O_{\sigma})=0$. (By the way, if $c^2=\lambda$ then $c^{-1}
\sigma(ct)$
is, clearly $\pm Cheb_d(t)$.)
\end{proof}

\begin{lem}
\label{sapte}
Assume $\sigma:\bP^1 \ra \bP^1$ has  degree $d \geq 2$ and admits a
non-zero semi-invariant form $\omega$ of positive 
weight not divisible by $p$. Assume either $p>d$ or $p=0$.
Assume moreover that $g(X_{\omega})=1$.
Then $\sharp P_{\sigma} \leq 4$ and
$\chi(O_{\sigma})=0$.
\end{lem}

\begin{proof}
Consider the lifting $\tau:X:=X_{\omega} \ra X$
  in Lemma \ref{patru} and let $G=\<\gamma\>$ ($n=\sharp G$)
 be the Galois group
of the covering $\pi:X \ra \bP^1$ as in the discussion following Equation
  \ref{gamma}. Then $\tau$ is \'{e}tale.
Let $S:=\pi(\Omega_{\pi})=\{A_1,...,A_s\}$ and 
$\tA_i \in \pi^{-1}(A_i)$. An easy analysis of the fixed points of
  $\gamma$ shows that the tuple
\[(e_{\pi}(\tA_1),...,e_{\pi}(\tA_s))\]
coincides, up to order, with one of the tuples
\begin{equation}
\label{mabuse}
(2,2,2,2),\ (3,3,3),\ (4,4,2),\ (6,3,2),
\end{equation}
according as $n$ is $2,3,4$ or $6$.
Note that for any $A \in \bP^1$ we have
\[\mu_{\sigma}(A)=lcm \{ \frac{e_{\pi}(\tA)}{e_{\pi}(\tB)}\ |\ 
\tB \in X,\ \tA=\tau^m(\tB),\ m \geq 1, \ \pi(\tA)=A\};\]
here we are using the fact that, for $B:=\pi(\tB)$, we have
\begin{equation}
e_{\sigma^m}(B) \cdot e_{\pi}(\tB)=e_{\pi}(\tA) \cdot e_{\tau^m}(\tB)=
e_{\pi}(\tA).
\end{equation}
By the above equation $e_{\pi}(\tA)=1$ imples $e_{\pi}(\tB)=1$ hence
$\mu_{\sigma}(A)=1$ for all $A \not\in S$. Now we claim that
\begin{equation}
\label{cuchi}
\mu_{\sigma}(A_i)=e_{\pi}(\tA_i)
\end{equation}
for all $i$. In view of the allowed values for $e_{\pi}(\tA_i)$ listed
in Equation \ref{mabuse}
 one sees that Equation \ref{cuchi} implies
$\sharp P_{\sigma} \leq 4$ and
$\chi(O_{\sigma})=0$ which will end our proof.

In order to check Equation \ref{cuchi} we need to show that for any
$i$ there exists $m \geq 1$ and $\tB \in X$ such that
$\tau^m(\tB)=\tA_i$, $e_{\pi}(\tB)=1$. This is, however, clear because
\[\sharp \ \tau^{-m}(\tA_i)=\sharp\ Ker\ \tau^m \longrightarrow
\infty\ \ \ \ as\ \ \ \ m \longrightarrow \infty,\]
hence the set $\tau^{-m}(\tA_i) \backslash \Omega_{\sigma}$ is
non-empty for $m$ sufficiently big.
\end{proof}

We are ready to prove our Theorem.

\begin{proof}
Assume $\bar{\sigma}_{\wp}$ admits an invariant form
 of non-zero  weight for infinitely many $\wp$'s.
By  Lemma \ref{trei}, for infinitely many $\wp$'s,
 $\bar{\sigma}_{\wp}$ admits a semi-invariant form
 of positive weight
not divisible by the residual characteristic of $\wp$.
By Lemmas \ref{sase} and \ref{sapte} 
\[\sharp(P_{\bsigma_{\wp}}) \leq 4,\ \ \chi(O_{\bsigma_{\wp}})=0\]
for infinitely many $\wp$'s. This  immediately implies
$\chi(O_{\sigma})=0$
hence,  by Proposition \ref{DHTHM}, $\sigma$ is flat.
\end{proof}

\bibliographystyle{amsplain}

\end{document}